\documentclass{ws-procs9x6}
\usepackage{amsfonts, amsmath, amstext,amssymb}
\usepackage{DraTex, AlDraTex}
\usepackage{array}

\def\phi{\varphi}

\begin{document}

\title{DEFINABILITY IN THE REAL UNIVERSE\,\footnote{Article based on an invited talk at Logic Colloquium 2009, University of Sofia, Bulgaria, August 5th, 2009.}}

\author{S. BARRY COOPER\,\footnote{Research supported by U.K.\ EPSRC grant No.\ GR
/S28730/01, and by a Royal Society International Collaborative 
Grant.}}
\address{Department of Pure Mathematics\break
University of Leeds,  Leeds LS2 9JT, U.K.}


\begin{abstract}
Logic has its origins in basic 
questions about the nature of the real world and how we describe it.  
This article seeks to bring out the physical and epistemological 
relevance of some of the more recent technical 
work in logic and computability theory. 
\vskip 1mm\noindent
``If you are receptive and humble, mathematics will lead you by the hand. 
Again and again, when I have been at a loss how to proceed, I have just 
had to wait until I have felt the mathematics lead me by the hand. It has led 
me along an unexpected path, a path where new vistas open up, a path 
leading to new territory, where one can set up a base of operations, from which one can survey 
the surroundings and plan future progress."
\vskip 0.5mm
\noindent {\em -- Paul Dirac, 27 November, 1975. 
In Paul A.\ M.\ Dirac Papers, 
Florida State University Libraries, Tallahassee, Florida, USA, No.\ 2/29/17.}
\end{abstract}



\bodymatter

\begin{quote}

\end{quote}

\section{Introduction}

Logic has an impressive history of addressing very basic 
questions about the nature of the world we live in. At the same time, it has 
clarified concepts and informal ideas about the world, and gone on to 
develop sophisticated technical frameworks within which these can 
be discussed. Much of this work is little known or understood by non-specialists, 
and the significance of it largely ignored. While notions such as set, proof and 
consistency have become part of our culture, other very natural abstractions 
such as that of definability are unfamiliar and disconcerting, even 
to working mathematicians.  The widespread interest in G\"odel's \cite{Godel_1931,Godel_1934} 
incompleteness results and their frequent application, often in questionable ways, 
shows both the potential for logicians to say something important about the world, while 
at the same time illustrating the limitations of what has been achieved so far. 
This article seeks to bring out the relevance of some of the more recent technical 
work in logic and computability theory. Basic questions addressed include: 
How do scientists represent and establish control over information about the 
universe? How does the universe itself exercise control over its own development? 
And  more feasibly: How can we reflect that control via our scientific and mathematical 
representations.

Definability --- what we can describe in terms of what we are given in a particular 
language --- is a key notion. As Hans Reichenbach (Hilary Putnam is perhaps his best-known 
student) found in the 1920s onwards, formalising  definability in the real world comes into its 
own when we need to clarify and better understand the content of a hard-to-grasp  
description of reality, such as Einstein's theory of general relativity. Reichenbach's 
seminal work\cite{reichenbach} on axiomatising relativity has become an ongoing project, carried forward 
today by Istvan Nemeti, Hajnal Andreka and their co-workers (see, for example, 
Andr\'eka,  Madar\'asz,  N\'emeti and  Sz\'ekely\cite{AMNS}). One can think 
of such work as paralleling the positive developments that models of computation 
enabled during the early days of computer science, bringing a surer grip on   
practical computation. But computability theory 
also gave an overview of what can be computed in principle, with corresponding 
technical developments apparently unrelated to applications.  
The real-world relevance of most of this theory remains conjectural. 

The capture of natural notions of describability and real-world robustness via 
the precisely formulated ones of definability and invariance also brings a 
corresponding development of theory, which can be applied in different 
mathematical contexts. Such an application does not just bring interesting 
theorems, which one just adds to the existing body of theory with 
conjectural relevance. It fills out the explanatory framework to a point 
where it can be better assessed for power and validity. And it is this  
which is further sketched out below. The basic ingredients are the 
notions of definability and invariance, and a mathematical context 
which best describes the scientific description of familiar causal structure. 

\section{Computability versus Descriptions}

In the modern world, {\em scientists} look for theories that enable predictions, 
and if possible predictions of a computational character. Everyone else 
lives with less constrained descriptions of what is happening, and is likely to happen. 
Albert Einstein\cite{einstein} might have expressed the view  in 1950 that:
\begin{quote}
When we 
say that we understand a group of natural phenomena, we 
mean that we have found a constructive theory which 
embraces them.
\end{quote}
But in everyday life people commonly use informal language to describe expectations 
of the real world from which constructive or computational content is not even attempted. 
And there is a definite mismatch between the scientist's drive to extend the reach of his or her 
methodology, and the widespread sense of an intrusion of algorithmic 
thinking into areas where it is not valid. A recent example is the controversy around 
Richard Dawkins' book\cite{dawkins06} on {\em The God Delusion}. This dichotomy has some basis in 
theorems from logic (such as G\"odel's incompleteness theorems): but the basis is more 
one for argument and confusion than 
anything more consensual. Things were not always so. 

If one goes back before the time of Isaac Newton, before the scientific era,  informal 
 {\em descriptions} of the nature of reality were 
the common currency of those trying to reason about the world. This might even impinge 
on mathematics --- as when the Pythagoreans wrestled with the ontology of irrational numbers. 
Calculation had a quite specific and limited role in society. 

\section{Turing's Model and \underline{In}computability}

In 1936, Turing\cite{Tu36} modelled what he understood of how a 
then human ``computer" (generally a young woman) might 
perform calculations ---  laying down rules that were very restrictive in a 
practical sense, but which enabled, as he plausibly argued, all that might be achieved 
with apparently more powerful computational actions. Just as the Turing machine's 
primitive actions (observing, moving, writing) were the key to modelling 
complex computations, so the Turing machine itself provided a route to the 
modelling of complex natural processes within structures 
which are discretely (or at least countably) presented. In this sense, 
it seemed we now had a way of making concrete the Laplacian model of 
science which had been with us in some form or other ever since the significance 
of what Newton had done become clear. 

But  the techniques for presenting a comprehensive range of computing machines gave us 
the {\em universal}  
Turing machine, so detaching computations from their material embodiments:  
and --- a more uncomfortable surprise --- by adding a quantifier to the perfectly down-to-earth description of the universal 
machine we get (and Turing\cite{Tu36} proved it) an \underline{in}computable object, the 
halting set of the machine. In retrospect, this becomes a vivid indication of how 
natural language has both an important real-world role, and quickly outstrips 
our computational reach. The need then becomes to track down 
material counterpart to 
 the simple mathematical schema which give rise to incomputability. 
Success provides a link to a rich body of theory and opens a Pandora's 
box of new perceptions about the failings of science and the 
nature of the real universe.

\section{The Real Universe as Discipline Problem}

The Laplacian model has a deeply ingrained hold on the rational mind. 
For a bromeliad-like late flowering of the paradigm we tend to think of Hilbert and 
his assertion of very general expectations for axiomatic mathematics. Or of the state of 
physics before quantum mechanics. The problem is that modelling the universe 
is definitely not an algorithmic process, and that is why intelligent, educated people can 
believe very different things, even in science. Even in mathematics.  So for many, 
the mathematical phase-transition from computability to incomputability,   
which a quantifier provides, is banned from the real 
world (see for example  Cotogno\cite{cot09}). However simple the {\em mathematical} route to incomputability, 
when looking out at the natural world, the trick is to hold the eye-glass to an unseeing eye. 
The global aspect of causality so familiar in mathematical structures is denied a connection 
with reality, in any shape or form. For a whole community, the discovery of 
incomputability made the real universe a real discipline problem. When Martin Davis\cite{davis04} says:
\begin{quote}
The great success of modern computers as all-purpose algorithm-executing engines embodying Turing's universal computer in physical form, makes it extremely plausible that the abstract theory of computability gives the correct answer to the question ÔWhat is a computation?Õ, and, by itself, makes the existence of any more general form of computation extremely doubtful.
\end{quote}
we have been in the habit of agreeing, in a mathematical setting. But in the 
context of a general examination of hypercomputational propositions (whatever 
the validity of the selected examples) it gives the definite impression of a   
defensive response to an uncompleted paradigm change. For convenience, we call this 
response\cite{davis06} --- that  ``there is no such discipline as hypercomputation"  -- 
``Davis' Thesis". 

The universal Turing machine freed us from the need actually {\em embody} 
the machines needed to host different computational tasks. The  
importance of this for building programmable computers was immediately 
recognised by  John 
von Neuman, and played a key role in the early history of the computer (see Davis\cite{davis:book}). 
The notion of 
a {\em virtual machine} is a logical extension of this tradition, which has found 
widespread favour amongst computer scientists and philosophers of a functionalist turn 
of mind --- for instance, there is the Sloman and Chrisley\cite{SC03} proposition 
for releasing consciousness from the philosophical inconvenience of embodiment (see also 
Torrance, Clowes and Chrisley\cite{TCC07}). Such attempts to tame nature are protected by a dominant 
paradigm, but there is plenty of dissatisfaction with them based on respect for the 
complex physicality of what we see. 

\section{A Dissenting Voice \dots }

Back in 1970, Georg Kreisel  considered one of the simplest physical 
situations presenting mathematical predictive problems. Contained within the 
mathematics one detects uncompleted infinities of the kind necessary for incomputability 
to have any significance for the real world. In a footnote to Kreisel\cite{Kr70} he proposed 
a collision problem related to the 
3-body problem, which might result in ``an analog computation of a 
non-recursive function". 

Even though Kreisel's view was built on many hours of deep thought 
about extensions of the Church-Turing thesis to the material universe --- much of this 
embodied in Odifreddi's 20-page discussion of the Church-Turing thesis in his 
book\cite{Od89} on Classical Recursion Theory --- it is not backed up by any proof of 
of the inadequacy of the Turing model 
 built on a precise description of the collision problem. 
 
 This failure 
 has become a familiar one, what has been described as a failure to 
 find `natural' examples of incomputability other than those computably equivalent 
 to the halting problem for a universal Turing machine --- with even that not 
 considered very natural by the mainstream mathematician. One requirement 
 of a `natural' incomputable set is that it be computably enumerable, like the set of 
 solutions of a diophantine equation, or the set of natural numbers $n$ such that 
 there exists a block of precisely $n$ 7s in the decimal expansion 
 of the real number $\pi$ --- or like the halting set of  a 
 given Turing machine. The problem is that 
 given a computably enumerable 
 set of numbers, there are essentially two ways of knowing its 
 incomputability. One way is to have designed the set oneself to have complement different to 
 any other set on a standard list of computably enumerable sets. Without 
 working relative to some other incomputable set, one just gets 
 canonical sets computably equivalent to the halting set of the universal Turing machine. 
 Otherwise the set one built has no known robustness, no definable 
 character one can recognise it by once it is built. 
 The other way of knowing a particular computably enumerable set to be incomputable 
 is to be able to compute one of the sets built via way one from the 
 given set. But only the canonical sets have been found so far to work in this way. 
 So it is known that there is a whole rich universe of computably inequivalent 
 computably enumerable sets --- but the only individual ones recognisably so 
 are computably equivalent to the halting problem. Kreisel's failure 
 is not so significant when one accepts that an arbitrary set picked from nature 
 in some way is very unlikely to be a mathematically canonical object. It seems 
 quite feasible that there is a mathematical theorem waiting to be proved, 
 explaining why there is no accessible procedure for verifying incomputability 
 in nature. 
 
 Since Kreisel's example, there have been other striking instances 
 of infinities in nature with the potential for hosting incomputability. 
 In {\em Off to Infinity  in Finite Time}  Donald Saari and Jeff Xia\cite{SX95}  describe 
 how one can even derive  singularities arising from the behaviour of 5 bodies 
 moving under the influence of the familiar Newtonian inverse square law. 
 
 There is a range of more complex examples which are 
 hard to fit into the standard Turing model, 
 ones with more real-world relevance. There is the persistence of problems 
 of predictability in a number of contexts. There is quantum uncertainty, constrained by 
 computable probabilities, but hosting what looks very much like randomness; there 
 are apparently emergent phenomena in many environments;  and chaotic 
 causal environments giving rise to strange attractors; 
 and one has  relativity and singularities (black holes), whose singular 
 aspects can host incomputability; .  Specially interesting is the 
renewed interest in analog and hybrid computing machines,  
  leading  Jan van Leeuwen and Jiri Wiedermann\cite{LW} to observe  
  that ``dots the classical Turing paradigm may no longer be fully appropriate to capture all features of present-day computing." And --- see later --- there is mentality, consciousness, 
  and the observed shortcomings of the mathematical models of these. 
 
 The disinterested observer of Martin Davis' efforts to keep nature contained within 
 the Turing/Laplacian model might keep in mind the well-known 
 comment of Arthur C.\ Clarke\cite{acc62} (Clarke's First Law)  that:
 \begin{quote}
 When a distinguished but elderly scientist states that something is possible, he is almost certainly right. When he states that something is impossible, he is very probably wrong.
 \end{quote}
 
 In what follows we look in more detail at three key challenges to the attachment 
 of Davis, and of a whole community, to 
  the Turing model in the form of Davis' thesis. 
  
  There is a reason for this. At first sight, it may seem unimportant to 
  know whether we have computational or predictive difficulties due to mere  
  complexity of a real-world computational task, or because of its actual incomputability. 
  And if there is no distinguishable difference between the two 
  possibilities, surely it cannot matter which pertains. Well, no. Attached to 
  two different mathematical characterisations one would expect different 
  mathematical theories. And there is a rich and well-developed theory of 
  incomputability. This mathematics may well constrain and give global 
  form to the real-world which it underlies. And these constraints and 
  structurings may be very significant for our experience and understanding 
  of the universe and our place in it. 

 \section{The Quantum Challenge}
 
 In the early days of quantum computing, there was some good news for 
 Davis' thesis from one of its most prominent supporters. David Deutsch, was one of the 
 originators of the standard model of quantum computation. In his  seminal 
 1985 article\cite{deutsch85}   
 ``Quantum Theory, the Church-Turing Principle and the Universal Quantum Computer" in the 
 {\em Proceedings of the Royal Society of London}, he introduced the notion of a `universal 
 quantum computer', and described how it might exploit quantum parallelism to compute 
 more efficiently than a classical Turing machine. But Deutsch is quite clear that 
 real computers based on this model would not compute anything not computable classically by 
 a Turing machine. And, of course, there are many other instances of successful reductions 
 of ``natural examples" of nature-based computational procedures to the Turing model. 
 
 But like Martin Davis, Deutsch\cite{deutsch06} is keen to take things further  ---  a { lot} further, 
 attempting a reduction of human mentality to the Turing model in a way even Turing in his 
 most constructive frame of mind might have had misgivings about:
 \begin{quote}
 I am sure we will have [conscious computers], I expect they will be 
 purely classical, and I expect that it will be a long time in the future. 
 Significant advances in our philosophical understanding of what consciousness is, 
 will be needed.
 \end{quote}
 
 Be this as it may, there are aspects of the underlying physics which are 
 not fully used in setting up the  standard model for quantum computing.  It is true that 
 {\em measurements} do play a role in a quantum computation, but in a tamed 
 guise. This is how Andrew Hodges explains it, in his 
 article {\em What would Alan Turing have done after 1954?} in the Teuscher volume\cite{teuscher}:
 \begin{quote}
 Von NeumannÕs axioms distinguished the {\bf U} (unitary evolution) and 
 {\bf R} (reduction) rules of quantum mechanics. Now, quantum computing
so far (in the work of Feynman, Deutsch, Shor, etc) is based on the {\bf U} process 
 and so computable. It has not made serious use of the {\bf R} process: 
 the unpredictable element that comes in with reduction, measurement, or 
 collapse of the wave function.
 \end{quote}
 The point being that measurements in the quantum context are  intrusive, with 
 outcomes governed by computable probabilities, but with 
 the mapping out of what goes on within those probabilities giving the 
 appearance of randomness. There are well-established formalisations 
 of the intuitive notion of randomness, largely coincident, and a large body of mathematical 
 theory built on these (see, for example,  Chaitin\cite{Ch87}, Downey and Hirschfeldt\cite{DH10}, 
 Nies\cite{Ni09}). A basic feature of the theory is the fact that randomness implies 
 incomputability (but not the converse). Calude and Svozil    \cite{CS08} have  
 extracted a suitable mathematical model of quantum randomness, built upon assumptions 
 generally acceptable to the physicists. Analysing the computability-theoretic 
 properties of the model, they are able to show that quantum randomness 
 does exhibit incomputability. But, interestingly, they are unable as yet to confirm that 
 quantum randomness 
 is mathematically random. 
 
 But quantum mechanics does not just present one of the toughest challenges to Davis' thesis.  
 It also presents the observer with a long-standing challenge to 
 its own  {\em realism}. Interpretations of the theory generally fail to satisfy everyone, 
 and the currently most widely accepted interpretations 
 contain what must be considered metaphysical assumptions. 
 When we have assembled the key ingredients, we will be in a position 
 to argue that the sort of fundamental thinking needed to 
 rescue the the theory from such assumptions is based on some very basic mathematics. 
 
  \section{Schr\"odinger's Lost States,  and the Many-Worlds Interpretation}
 
 One way of describing the quantum world is via the Schr\"odinger wave equation. 
 What Hodges refers to above are the processes for change of the wave equation describing 
 the quantum state of a physical system. On the one hand, one has 
 deterministic continuous evolution via Schr\"odingerÕs equation,  involving 
 superpositions of basis states. On the other, one has 
probabilistic non-local discontinuous change due to measurement. 
 With this, one observes a jump to a single basis state. The interpretive question 
 then is:  {\em Where do the other states go?}
 
 Writing with hindsight: If the physicists knew enough logic, they would have been able to make a good guess. 
 And if the logicians had been focused enough on the foundations of  quantum mechanics 
 they might have been able to tell them. 
 
 As it is,  physics 
 became a little weirder around 1956. The backdrop to this is the sad and 
 strange life-story of Hugh Everett III and his family, through which strode the 
 formidable 
 John Wheeler, Everett's final thesis advisor, and Bryce DeWitt , 
 who in 1970 coined the term `Many-Worlds'  for Everett's 
 neglected and belittled idea: an idea whose day came too late to help   
 the Everett family, now only survived by the son Mark who relives 
 parts of the tragic story 
 via an autobiography \cite{Ev08} and 
 appropriately left field confessional creations as leader of  the Eels rock band. 
 
 Many-Worlds, with a little reworking, did away with the need to explain the transition from 
 many superposed quantum states to the `quasi-classical' uniqueness 
 we see around us. The multiplicity survives and permeates micro-  to macro-reality, 
 via a 
 decohering bushy branching of alternative histories, with us 
 relegated to to our own self-contained branch. Max Tegmark has 
 organised the multiplying   
 variations on the Many-Worlds theme  into hierarchical levels of 
 `multiverses', from modest to more radical proposals, with even 
 the underlying mathematics and the consequent 
 laws of physics individuating at Level IV. 
 Of course, if one does not 
 bother anymore to explain why our universe works so 
 interestingly, one needs 
 the `anthropic principle' on which to base our experience of the world --- 
 ``We're here because we're here because we're here because we're here \dots ", 
 as they sang during the Great War, marching towards the trenches. 
 The attraction of this picture derives from the drive for 
 a coherent overview, and the lack of a better one. 
 As David Deutsch put it in {\em The Fabric of Reality} \cite[p.48]{deutsch97}:
 \begin{quote}
 \dots  understanding the multiverse is a precondition for 
 understanding reality as best we can. Nor is this said in a spirit of 
 grim determination to seek the truth no matter how unpalatable it may 
 be \dots It is, on the contrary, because the resulting world-view is so 
 much more integrated, and makes more sense in 
 so many ways, than any previous world-view, and 
 certainly more than the cynical pragmatism which 
 too often nowadays serves as surrogate for a world-view amongst scientists.
 \end{quote}
 
 Here is a very different view of the multiverse from the distinguished 
 South African mathematician George Ellis \cite[p.198]{ellis96}, one-time collaborator 
 of  Stephen Hawking:
 \begin{quote}
 The issue of what is to be regarded as an ensemble of `all possible'  
 universes is unclear, it can be manipulated to produce any result you want \dots The argument that this infinite ensemble actually exists can be claimed to have a certain explanatory economy (Tegmark 1993), although others would claim that Occam's razor has been completely abandoned in favour of a profligate excess of existential multiplicity, extravagantly hypothesized in order to explain the one universe that we do know exists.
 \end{quote}
 
 The way out of this foundational crisis, as with previous ones in mathematics and science, 
 is to adopt a more constructive approach. In this way, one can combine the 
 attractions of Tegmark's \cite{Teg08} {\em Mathematical Universe Hypothesis} (MUH) with the 
 discipline one gets from the mathematics of what can be built from very small beginnings.

  \section{Back in the One World \dots }
 
 A constructive approach is not only a key to clarifying the 
  interpretive problem. Eliminating the redundancy 
  of parallel universes, and the reliance on the anthropic principle, 
  also entails the tackling of the unsatisfactory 
  arbitrariness of various aspects of the standard model. 
The exact values of the constants of nature, 
subatomic structure, the geometry of space --- 
all confront the standard model of particle physics 
with a foundational problem.   Alan Guth, inventor 
of the `cosmic inflation'  needed to make sense of our picture of  
the early universe, asks \cite{guth97}:
\begin{quote}
If the creation of the universe can be described as a quantum process,
we would be left with one deep mystery of existence: What is it that determined 
the laws of physics?
\end{quote}
 And Peter Woit, in his recent book   \cite{woit06} 
{\em  Not Even Wrong --- The Failure of String Theory and the 
Continuing Challenge to Unify the Laws of Physics}, comments on the arbitrary constants 
one needs to give the right values to get the standard model to behave properly:
\begin{quote}
One way of thinking about what is unsatisfactory about the standard model is that it leaves seventeen non-trivial numbers still to be explained, \dots 
\end{quote}
 Even though the exact number of constants undetermined by  theory, but needing special 
 fine-tuning to make the standard model fit with observation, does vary,  
even one is too many.
  This dissatisfaction with aspects of the standard model goes back to Einstein. Quoting 
 from Einstein's {\em Autobiographical Notes} \cite[p.63]{Ei69}:
 \begin{quote}
 \dots I would like to state a theorem which at present 
 can not be based upon anything more than upon a faith in the 
 simplicity, i.e. intelligibility, of nature \dots nature is so constituted 
 that it is possible logically to lay down such strongly determined 
 laws that within these laws only rationally completely 
 determined constants occur (not constants, therefore, whose numerical 
 value could be changed without destroying the theory) \dots
 \end{quote}
 
 What is needed is mathematics which does more than express mechanistic 
 relationships between basic entities. One needs theory expressed in 
 language strong enough to encapsulate not just relations on the material world,  
 but relations on such relations ---  relations which entail qualifications 
 sophisticated enough to determine all aspects of the our universe, including 
 the laws of nature themselves. 
 Or, as Roger Penrose terms it \cite[pp.106-107]{Pe87}, we need to capture {\em Strong Determinism}, 
 whereby : 
 \begin{quote}
  \dots all the complication, 
 variety and apparent randomness that we see all about 
 us, as well as the precise physical laws, are all exact and 
 unambiguous consequences of one single coherent mathematical structure.
 \end{quote}
 The article \cite{CCSS} of Calude, Campbell, Svozil and Stefanescu on 
 {\em Strong determinism vs. computability} contains a useful discussion of the 
 computability-theoretic ramifications of strong determinism. 
 
 In the next section we examine some 
  more approachable phenomena than those at the quantum level. 
  Even though the challenge these present 
  to Davis' Thesis is less obvious than that of quantum uncertainty, they do  
    point us in the direction of the 
 mathematics needed to make sense of strong determinism. 
 
 \section{The Challenge from Emergence}
 
 The waves on the seashore, the clouds scudding across the sky,  the 
 complexity of the Mandelbrot set ---   observing these, 
 one is made aware of limits on what we can practically compute. The underlying 
 rules governing them are known, but that is not enough. When we talk about the problem of 
 `seeing the wood for the trees'  we are approaching the 
  gap between micro- and macro- events from another direction.  Either way, 
  there are commonly encountered situations in which either reduction, or 
  seeing the `big picture', entails more than a computation. 
  
  Although an interest in such things goes back to Poincar\'e --- we already 
  mentioned the 3-body problem --- it was the  second 
  half of the twentieth 
  century saw the growth of chaos theory, and a greater of awareness of the 
  generation of informational complexity via simple rules, accompanied by the 
  {\em emergence} of new regularities. The most mundane and apparently 
  uncomplicated situations could provide examples --- such as Robert Shaw's \cite{shaw84}  
  strange attractor arising from an appropriately paced dripping tap. 
  And inhospitable as turbulent fluids might appear, there too 
  higher order formations might emerge and be subject to 
  mathematical description, as demonstrated by David Ruelle 
  (see Ruelle\cite{ruelle93}) another 
  early pioneer in the area. Schematic metaphors for such 
  examples are provided by the cellular 
  automaton (CA) model, and famously by John Conway's Game of Life. Here is the 
  musician Brian Eno \cite{eno} talking in relation to how his 
   creative work on `generative music' was influenced by `Life':
  \begin{quote}
  These are terribly simple rules and you would think it probably couldn't produce anything very interesting. Conway spent apparently about a year finessing these simple rules. 
  \dots  He found that those were all the rules you needed to produce something that appeared life-like. \par
  What I have over here, if you can now go to this Mac computer, please. I have a little group of live squares up there. When I hit go I hope they are going to start behaving according to those rules. There they go. I'm sure a lot of you have seen this before. What's interesting about this is that so much happens. The rules are very, very simple, but this little population here will reconfigure itself, form beautiful patterns, collapse, open up again, do all sorts of things. It will have little pieces that wander around, like this one over here. Little things that never stop blinking, like these ones. What is very interesting is that this is extremely sensitive to the conditions in which you started. If I had drawn it one dot different it would have had a totally different history. This is I think counter-intuitive. One's intuition doesn't lead you to believe that something like this would happen.
  \end{quote}
  Margaret Boden and Ernest Edmonds \cite{BE09} make a case for  
 generative art,  emergent from automata-like 
 computer environments, really qualifying as art. 
  While computer pioneer Konrad Zuse was impressed 
  enough by the potentialities of cellular automata to 
   suggest \cite{zuse} 
  that 
  the physics of the universe might be CA computable. 
  
  A specially useful key to a general mathematical understanding of such phenomena 
  is the well-known link between emergent structures in 
  nature, and familiar mathematical objects, such as the Mandelbrot and Julia sets. 
  These mathematical metaphors for real-world complexity 
  and associated patterns have caught the attention of many --- such as  
  Stephen Smale \cite{BCSS} and Roger Penrose --- 
  as a way of getting a better grip on the computability/complexity of emergent 
  phenomena. Here is Penrose \cite{penrose94} describing his fascination with the Mandelbrot 
  set:
  \begin{quote}
  Now we witnessed \dots a certain extraordinarily complicated 
  looking set, namely the Mandelbrot set. Although the rules 
  which provide its definition are surprisingly simple, 
  the set itself exhibits an endless variety of highly elaborate structures.
   \end{quote}
   
  As a mathematical analogue of emergence in nature, what are the  
  distinctive mathematical characteristics of the Mandelbrot set? 
  It is derived from a simple polynomial formula over the 
  complex numbers, via the addition of a couple of quantifiers. 
  In fact, with a little extra work, the quantifiers can be reduced to just one. 
  This gives the definition the aspect of a familiar object from classical 
  computability theory --- namely, a $\Pi^0_1$ set. Which is just the 
  level at which we might not be surprised to encounter incomputability. 
  But we have the added complication of working with real (via complex) numbers 
  rather than just the natural numbers. This creates room for a certain amount 
  of controversy around the use of the BSS model of real computation 
  (see Blum, Cucker, Shub and Smale\cite{BCSS}) to show the incomputability of the Mandelbrot set and most 
  Julia sets. The 2009 book by Mark Braverman and Michael Yampolsky \cite{BY09} on 
  {\em Computability of Julia Sets} is a reliable guide to recent 
  results in the area, including those using the more mainstream computable analysis 
  model of real computation. The situation is not 
  simple, and the computability of the Mandelbrot set, as of 
  now,  is still an open question.
  
  What is useful, in this context, is that these examples both connect with 
  emergence in nature, and share logical form with well-known 
  objects which transcend the standard Turing model. As such, they point to 
  the role of extended language in a real context taking us beyond 
  models which are purely mechanistic. And hence give us a route to 
  mathematically capturing the origins of emergence in nature, and to 
  extending our understanding of how nature computes. We can now view 
  the halting set of a universal Turing machine as an emergent phenomenon, 
  despite it not being as pretty visually as our Mandelbrot and Julia examples. 
 
 One might object that there is no evidence that 
 quantifiers and other globally defined operations have any 
 existence in nature beyond the minds 
 of logicians. But how does nature know anything about any logical construct? 
 The basic logical operations derive their basic status from their association 
 with elementary algorithmic relationships over information. Conjunction 
 signifies an appropriate and very simple merging of information, of the kind 
 commonly occurring in nature. Existential quantification expresses 
 projection, analogous to a natural object throwing a shadow on a bright 
 sunny day. And if a determined supporter of Davis' Thesis plays at god, and 
 isolates a computational environment with the aim of bringing it 
 within the Turing model, then the 
 result is the delivery of an identity to that environment, the creating a 
 natural entity --- like a human being, perhaps ---  with undeniable 
 naturally emergent global attributes. 
 
 There are earlier, less schematic approaches to the mathematics of emergence. 
 Ones which fit well with the picture so far. 
 
  It often happens that  when one gets interested in a particular aspect of computability, 
 one finds Alan Turing was there before us. Back in the 1950s, Turing 
 proposed a simple reaction-diffusion system describing chemical reactions and 
 diffusion to account for morphogenesis, i.e., the 
 development of form and shape in biological systems. One can find a full 
 account of the background to Turing's seminal intervention in the field 
 at Jonathan Swinton's well-documented webpage \cite{swinton} on  {\em 
 Alan Turing and morphogenesis}. One of Turing's main achievements was to 
 come up with mathematical descriptions --- differential equations --- 
 governing such phenomena 
 as Fibonacci phyllotaxis: the surprising showing of 
 Fibonacci progressions in such things as the criss-crossing spirals 
 of a sunflower head. As Jonathan Swinton describes:
 \begin{quote}
 In his reaction-diffusion system [Turing] had the first and one of the 
 most compelling models mathematical biology has devised 
 for the creation process. In his formulation of the Hypothesis of 
 Geometrical Phyllotaxis he expressed simple rules adequate for 
 the appearance of Fibonacci pattern. In his last, unfinished work he 
 was searching for plausible reasons why those rules might hold, and it 
 seems only in this that he did not succeed. It would take many decades 
 before others, unaware of his full progress, would retrace his steps 
 and finally pass them in pursuit of a rather beautiful theory.
 \end{quote}
 Most of Turing's work in this area was unpublished in his lifetime, 
  only appearing in 1992 in the {\em Collected Works} \cite{Turing:1992}. 
  Later work, coming to fruition just after Turing died, was carried forward by his student Bernard 
  Richards, appearing in the thesis \cite{Ri54}. See Richards\cite{richards} for an 
  interesting account of Richards' time working with Turing.
  
  The field of {\em synergetics}, 
  founded by the German physicist Hermann Haken, provides another  mathematical approach to emergence. Synergetics is a multi-disciplinary approach to 
 the study of the origins and evolution of macroscopic patterns and spacio-temporal structures in interactive systems. 
   An important feature of synergetics for our purposes 
  is its focus on  {\em self-organizational} processes in science and the humanities, particularly 
  that of autopoiesis.  An instance of an 
  autopoietic system is a biological cell, and is distinguished by being   
  sufficiently autonomous and operationally closed, to recognisably self-reproduce.   
  
  A particularly celebrated example of the technical effectiveness of the theory  is 
  Ilya Prigogine's achievement of the Nobel Prize for Chemistry in 1977 for his 
  development of dissipative structure theory and its 
  application to thermodynamic systems far from equilibrium, with subsequent 
  consequences for self-organising systems. Nonlinearity and irreversibility are 
  associated  key aspects of the processes modelled  
  in this context. 
  
  See Michael Bushev's comprehensive review of the 
  field in his book \cite{bushev} {\em Synergetics -- Chaos, Order, Self-Organization}. 
  Klaus Mainzer's book \cite{mai94} on {\em Thinking in Complexity: The 
  Computational Dynamics of Matter, Mind, and Mankind} puts synergetics 
  in a wider context, 
  and mentions such things as synergetic computers. 
   
  The emphasis of the synergetists on {\em self-organisation} in relation to the emergence of order from 
  chaos is important in switching attention from the {\em surprise} highlighted by so many 
  accounts of emergence, to the {\em autonomy} and {\em internal organisation} 
  intrinsic to the phenomenon. People like Prigogine found within synergetics, 
  as did Turing for morphogenesis, precise descriptions of  previously 
  mysteriously emergent order. 
  
 \section{A Test for Emergence}
 
 There is a problem with the big claims made for emergence in many different 
 contexts. Which is that, like with `life', nobody has a good definition of it. Sometimes, 
 this does matter. Apart from which history is littered with instances of vague concepts clarified 
 by science, with huge benefits to our understanding of the world and to the progress of 
 science and technology. The clarification of what we mean by a 
 {\em computation}, and the subsequent development of the computer 
 and computer science is a specially relevant example here. Ronald C. Arkin, in his book 
 \cite[p.105]{arkin} {\em Behaviour-Based Robotics}, summarises the problem as 
 it relates to emergence:
 \begin{quote}
  Emergence is often invoked in an almost mystical 
  sense regarding the capabilities of behavior-based systems. 
  Emergent behavior implies a holistic capability where the sum is 
  considerably greater than its parts. It is true that what occurs in a 
  behavior-based system is often a surprise to the system's designer, 
  but does the surprise come because of a shortcoming 
  of the analysis of the constituent behavioral building blocks and their 
  coordination, or because of something else?
 \end{quote}
 
 There is a salutary warning from the history of British Emergentists, who had their heyday   
  in the early 1920s --- Brian McLaughlin's 
   book \cite{mclaughlin}. The notion of emergence has been found to be a useful concept 
 from at least the time of John Stuart Mill, back in the nineteenth century. 
   The emergentists of the 1920s used the concept to explain the 
   irreducibility of the `special sciences', postulating a hierarchy with 
   physics at the bottom, followed by chemistry, biology, social science etc. 
   The emergence was seen, anticipating modern thinking, as being irreversible, imposing the irreducibility 
   of say biology to quantum theory. Of course  the British emergentists experienced  
 their heyday before the great quantum discoveries of the late 1920s, and as 
 described in McLaughlin\cite{mclaughlin}, this was in a sense their undoing. One of the leading figures 
 of the movement was the Cambridge philosopher C.\ D.\ Broad, described by Graham Farmelo 
 in his biography of Paul Dirac \cite[p.39]{farmelo} as being in 1920 
 ``one of the most talented young philosophers 
 working in Britain". In many ways a precursor of the current philosophers 
 arguing for the explanatory role of emergence in the philosophy of mind, 
 Charlie Broad was alive to the latest scientific developments, lecturing to 
 the young Paul Dirac on Einstein's new theory of relativity while they were both at Bristol. 
 But here is Broad writing in 1925 \cite[p.59]{broad} about the `emergence' of salt crystals:
 \begin{quote}
 \dots the characteristic behaviour of the whole \dots could not, even in 
 theory, be deduced from the most complete knowledge of the 
 behaviour of its components \dots This \dots is what I understand 
 by the `Theory of Emergence'. I cannot give a conclusive example of it, 
 since it is a matter of controversy whether it actually applies to anything \dots 
 I will merely remark that, so far as I know at present, the characteristic 
 behaviour of Common Salt cannot be deduced from the most 
 complete knowledge of the properties of Sodium in isolation; or of 
 Chlorine in isolation; or of other compounds of Sodium, \dots 
 \end{quote}
 The date 1925 is significant of course, it was in the years following 
 that Dirac and others developed the quantum mechanics which would 
 explain much of chemistry in terms of locally described interactions 
 between sub-atomic particles. The reputation of the emergentists, for 
 whom such examples had been basic to their argument for the 
 far-reaching relevance of emergence, never 
 quite recovered. 
 
 For Ronald, Sipper and Capcarr\`ere in 1999, Turing's approach to 
 pinning down intelligence in machines suggested a test for 
 emergence. Part of the thinking would have been that emergence, 
 like intelligence, is something we as observers think we can recognise;  
 while the complexity of what we are looking for resists 
 observer-independent analysis. The lesson is to police the 
 observer's evaluation process, laying down some optimal rules 
 for a human observer. Of course, the Turing Test is specially 
 appropriate to its task, our own experience of human 
 intelligence making us well-qualified to evaluate the putative 
 machine version. Anyway, the Emergence Test of 
 Ronald, Sipper and Capcarr\`ere \cite{RSC} for emergence 
 being present in a system, modelled on  
 the Turing Test, had the following three ingredients:
 \begin{enumerate}
 
\item {\bf Design:} The system has been constructed by the designer, 
by describing local elementary interactions between components 
(e.g., artificial creatures and elements of the environment) in a 
language ${\frak L}_{1}$.

\item {\bf Observation:} The observer is fully aware of the design, 
but describes global behaviors and properties of the 
running system, over a period of time, 
using a language ${\frak L}_2$.

\item {\bf Surprise:} The language of design ${\frak L}_1$ and the 
language of observation ${\frak L}_2$ are distinct, and the 
causal link between the elementary interactions programmed 
in ${\frak L}_1$ and the behaviors observed in ${\frak L}_2$ is non-obvious to the 
observer --- who therefore experiences surprise. In other 
words, there is a cognitive dissonance between the observer's 
mental image of the system's design stated in ${\frak L}_1$ and his 
contemporaneous observation of the system's behavior stated in 
${\frak L}_2$.

\end{enumerate}

Much of what we have here is what one would expect, extracting 
the basic elements of the previous discussion, and expressing 
it from the point of view of the assumed observer. But an 
 ingredient which should be noted is the formal distinction 
 between the language ${\frak L}_1$ of the design and that 
 of the observer, namely ${\frak L}_2$. This fits in with our 
 earlier mathematical examples: the halting set of a universal 
 Turing machine, and the Mandelbrot set, where the new language 
 is got by adding a quantifier --- far from a minor augmentation of 
 the language, as any logician knows. And it points to the 
 importance of the language used to describe the phenomena, 
 an emphasis underying the next section.

 \section{Definability the Key Concept}
 
We have noticed that it is often possible to get descriptions of emergent 
properties in terms of the elementary actions from which they arise. 
For example, this is what Turing did for the role of Fibonacci numbers 
in relation to the sunflower etc. This is not unexpected, it is characteristic 
of what science does. And  
in mathematics, it is well-known that complicated 
descriptions may take us beyond what is computable. 
This could be seen as a potential source of surprise in emergence.

But one can turn this viewpoint around, and get something 
more basic. There is an intuition that entities do not just 
generate descriptions of the rules governing them:
 they actually {\em exist} because of, and according to mathematical 
 laws. And that for entities that we can be aware of, 
 these will be mathematical laws which are susceptible to description. 
  That it is the describability that is key to their observability. But 
  that the existence of such descriptions is not enough to ensure we 
  can access them, even though they have algorithmic content 
  which provides the stuff of observation.
  
  It is hard to for one to say anything new. In this case Leibniz was there 
  before us, essentially with his Principle of Sufficient Reason. 
  According to Leibniz \cite{leibniz} in 1714:
  \begin{quote}
  \dots there can be found no fact that is true or existent, 
  or any true proposition, without there being a sufficient 
  reason for its being so and not otherwise, although we cannot 
  know these reasons in most cases.
  \end{quote} 
 Taking this a little further --  natural phenomena not only generate 
 descriptions, but arise and derive form from them.   
  And this  connects with a useful abstraction --  that of mathematical definability,
   or, more generally, invariance under the automorphisms of the appropriate structure. 
So giving precision to our experience of emergence as a potentially non-algorithmic determinant of events. 
   
  This is a familiar idea in the mathematical context. The relevance 
  of definability for the real world is implicitly present in Hans Reichenbach's 
  work \cite{reichenbach} on the axiomatisation of relativity. It was, of course, 
  Alfred Tarski who gave a precise logical form to the notion of definability. 
  Since then logicians have worked within many different mathematical 
  structures, succeeding in showing that different operations and relations 
  are non-trivially definable, or in some cases undefinable, in terms 
  of given features of the structure. Another familiar feature of mathematical 
  structures is the relationship between definability within the structure 
  and the decidability of its theory (see Marker \cite{marker}), giving 
  substance to the intuition that knowledge of the world is so hard to 
  capture, because so much can be observed and described. Tarski's proof of 
  decidability of the real numbers, contrasting 
  with the undecidability of arithmetic,  fits with the fact that 
  one cannot even define the integers in the structure of the real numbers. 
  
  Unfortunately, outside of logic, and certainly outside of mathematics, 
  the usefulness of definability remains little understood. And 
  the idea that features of the real world may actually be undefinable 
  is, like that of incomputability, a recent and unassimilated 
  addition to our way of looking at things.
  
  At times, definability or its breakdown  comes disguised within quite familiar 
  phenomena. In science, particularly in basic physics, 
  symmetries play an important role. One might be surprised at this, 
  wondering where all these often beautiful and surprising 
  symmetries come from. Maybe designed by some higher power? 
  In the context of a mathematics in which undefinability 
  and nontrivial automorphisms of mathematical structures is 
  a common feature, such symmetries  lose 
  their unexpectedness. When Murray Gell-Mann 
  demonstrated the relevance of SU(3) group symmetries 
  to the quark model for classifying of elementary particles, 
  it was based on lapses  
  in definability of the strong nuclear force
   in relation to quarks of differing  flavour. The automorphisms 
   of which such symmetries are an expression give a clear 
   route from fundamental mathematical structures and their 
   automorphism groups  to far-reaching 
   macro-symmetries in nature. If one accepts that such basic 
   attributes as {\em position} can be subject to failures of 
   definability, one is close to restoring realism to 
   various basic subatomic phenomena.  
   
   One further observation: Identifying emergent phenomena 
   with material expressions of definable relations suggests 
   an accompanying {\em robustness} of such phenomena. 
   One would expect the mathematical characterisation to 
   strip away much of the mystery which has made emergence so 
   attractive to theologically inclined philosophers of mind,  
   such as Samuel Alexander  \cite[p.14]{alexander}:
\begin{quote}
The argument is that mind has certain specific characters to which there is or even 
can be no neural counterpart \dots Mind is, according to our interpretation of the 
facts, an `emergent'  from life, and life an emergent from a lower physico-chemical 
level of existence.
\end{quote}
And further \cite[p.428]{alexander}:
\begin{quote}
In the hierarchy of qualities the next higher quality to the highest attained is deity. God is the whole universe engaged in process towards the emergence of this new quality, and religion is the sentiment in us that we are drawn towards him, and caught in the movement of the world to a higher level of existence.
\end{quote}

In contrast, here is Martin Nowak, 
Director of the Program for Evolutionary Dynamics at Harvard University, writing in the 
collection \cite{brockman} {\em What We Believe But Cannot Prove}, describing  
the sort of robustness we would expect:
\begin{quote}
I believe the following aspects of evolution to be true, without
knowing how to turn them into (respectable) research topics.\par \noindent
Important steps in evolution are robust. Multicellularity evolved at least ten times. There are several independent
origins of eusociality. There were a number of lineages leading from primates to humans. If our ancestors had not evolved language, somebody else would have.
\end{quote}
   
  What is meant by robustness here is that there is mathematical content 
  which enables the process to be captured and moved between different 
  platforms. Though it says nothing about the relevance of embodiment or 
  the viability of virtual machines hostable by canonical machines. 
  We return to this later. On the other hand, it gives us a handle 
  on representability of emergent phenomena, a key aspect 
  of intelligent computation.

 \section{The Challenge of Modelling Mentality}
 
 Probably the toughest environment in which to road-test the 
 general mathematical framework we have associated with 
 emergence is that of human mental 
 activity. What about the surprise ingredient of the 
 Emergence Test? 
 
 Mathematical thinking provides an environment in which major   
 ingredients -- Turing called them intuition and ingenuity, others might call them 
 creativity and reason -- are easier to clearly separate. 
 A classical source of information and analysis of such thinking is 
 the French mathematician 
 Jacques Hadamard's  
  {\em The Psychology of Invention in the Mathematical Field} \cite{hadamard}, 
  based on personal accounts   
  supplied by distinguished informants such as Poincar\'e, Einstein 
  and Polya. Hadamard was particularly struck by Poincar\'e's thinking, 
  including a 1908 address of his to the 
  French Psychological Society in Paris on the topic of 
   {\em Mathematical Creation}. Hadamard followed Poincar\'e and Einstein 
   in giving an important role to unconscious thought processes, and 
   their independence of the role of language and mechanical reasoning.  
   This is Hadamard's account, built 
   on that of  Poincar\'e\cite{poincare}, of Poincar\'e's 
   experience of struggling with 
   a problem:
   \begin{quote}
   At first Poincar\'e attacked [a problem] vainly for a fortnight, attempting to prove there could not be any such function \dots  [quoting Poincar\'e]:\par
``Having reached Coutances, we entered an omnibus to go some place or other. At the moment when I put my foot on the step, the idea came to me, without anything in my former thoughts seeming to have paved the way for it \dots I did not verify the idea 
   \dots I went on with a conversation already commenced, but I felt a perfect certainty. 
   On my return to Caen, for conscience sake, I verified the result at my leisure."
   \end{quote}
  This experience will be familiar to most research mathematicians -- the 
  period of incubation, the failure of systematic reasoning, and the 
  surprise element in the final discovery of the solution: a surprise that 
  may, over a lifetime, lose some of its bite with repetition and familiarity, but which 
  one is still compelled to recognise as being mysterious and worthy of surprise. 
  Anyway, the important third part of the Emergence Test is satisfied here.
  
  Perhaps even more striking is the fact that Poincar\'e's solution had that 
  robustness we looked for earlier: the solution came packaged and 
  mentally represented in a form which enabled it to be carried home and 
  unpacked intact when back home. Poincar\'e just carried on with his 
  conversation on the bus, his friend presumably unaware of the 
  remarkable thoughts coursing through the mathematicians mind. 
  
  Another such incident emphasises the lack of uniqueness and the 
  special character of such incidents -- Jacques Hadamard \cite{hadamard} 
  quoting Poincar\'e 
  again:
  \begin{quote}
  ``Then I turned my attention to the study of some arithmetical questions 
  apparently without much success \dots Disgusted with my failure, 
  I went to spend a few days at the seaside and thought of something 
  else. One morning, walking on the bluff, the idea came to me, with
   just the same characteristics of brevity, suddenness and immediate 
   certainty, that the arithmetic transformations of indefinite ternary
    quadratic forms were identical with those of non-Euclidian geometry."
  \end{quote}
 
 What about the {\em design}, and the observer's {\em awareness} of the design? 
Here we have a large body of work , most notably  from neuro-scientists and 
philosophers, and an increasingly detailed knowledge of the workings 
of the brain. What remains in question -- even  accepting the brain as the 
design (not as simple as we would like!) -- is the exact nature of the connection 
between the design and the emergent level of mental activity. This is an 
area where the philosophers pay an important role in clarifying 
problems and solutions, while working through consequences and 
consistencies. 

The key notion, providing a kind of workspace for working through 
alternatives, is that of {\em supervenience}. According to Jaegwon 
Kim \cite[pp.14--15]{kim98}, supervenience:
\begin{quote}
\dots represents the idea that mentality is at bottom 
physically based, and that there is no free-floating mentality 
unanchored in the physical nature of objects and events in which 
it is manifested.
\end{quote}
There are various formulations. This one is from the online {\em Stanford 
Encyclopedia of Philosophy}:
\begin{quote}
A set of properties {\bf A} supervenes upon another set 
{\bf B} just in case no two things can differ with respect to 
{\bf A}-properties without also differing
with respect to their {\bf B}-properties.
\end{quote}
So in this context, it is the mental properties which are thought to 
supervene on the neuro-physical properties. All we need to know is 
are the details of how this supervenience takes place. And what 
throws up difficulties is our own intimate experience of the 
outcomes of this supervenience. 

One of the main problems relating to supervenience is the so-called 
`problem of mental causation', the old problem which undermined  
the Cartesian conception of mind-body dualism. The persistent 
question is: {\em How can mentality have a causal role in a world that 
is fundamentally physical?} Another unavoidable problem is that 
of  `overdetermination' -- the problem of phenomena 
having both mental and physical causes. For a pithy expression of the 
problem, here is Kim \cite{kim05} again:
\begin{quote}
\dots the problem of mental causation is solvable only if mentality 
is physically reducible; however, phenomenal consciousness resists physical reduction, 
putting its causal efficacy in peril.
\end{quote}

It is not possible here, and not even useful, to go into the intricacies 
of the philosophical debates which rage on. But it is important 
to take on board the lesson that a crude 
mechanical connection between mental activity and 
the workings of the brain will not do the job. Mathematical 
modelling is needed to clarify the mess, but has to meet 
very tough demands. 
 
 \section{Connectionist Models to the Rescue?}
 
 Synaptic interactions are basic to the workings of the brain, and 
 connectionist models based on these are the first hope.  
 And there is optimism about such models from such leading figures in the field 
 as Paul Smolensky \cite{smol88}, recipient of the 
 2005 David E.\ Rumelhart Prize:
 \begin{quote}
 There is a reasonable chance that connectionist models will 
 lead to the development of
new somewhat-general-purpose self-programming, massively
  parallel analog computers, and a new theory of analog parallel 
  computation: they may possibly even challenge the strong 
  construal of Church's Thesis as the claim that the class of 
  well-defined computations is exhausted by those of Turing machines.
 \end{quote}
 And it is true that connectionist models have come 
 a long way since Turing's 1948   discussion \cite{tur48} 
 of `unorganised machines', and McCulloch and Pitts'  1943  
 early paper \cite{MP43} on neural nets. (Once again, 
 Turing was there at the beginning, see Teuscher's book \cite{teusch02} 
 on {\em Turing's Connectionism}.)
 
 But is that all there is? For Steven Pinker 
 \cite{pinker97}  ``\dots neural networks alone 
 cannot do the job". And focusing on our elusive higher functionality, 
 and the way in which mental images are recycled and incorporated 
 in new mental processes, 
 he points to a ``kind of mental fecundity called recursion":
 \begin{quote}
 We humans can take an entire proposition and give it a role in some 
 larger proposition. Then we can take the larger proposition and embed 
 it in a still-larger one. Not only did the baby eat the slug, but the 
 father saw the baby eat the slug, and I wonder whether the father 
 saw the baby eat the slug, the father knows that I wonder whether 
 he saw the baby eat the slug, and I can guess that the father 
 knows that I wonder whether he saw the baby eat the slug, and so on.
 \end{quote}
 
 Is this really something new? Neural nets can handle recursions 
 of various kinds. They can exhibit imaging and representational 
 capabilities. They can learn. The problem seems to be with modelling the 
 holistic aspects of brain functionalism. It is hard to 
 envisage a model at the level of neural networks which successfully 
 represent and communicate its own global informational structures. 
 Neural nets do have many of the basic ingredients of what one observes 
 in brain functionality, but the level of developed synergy of the ingredients 
 one finds in the brain does seem to occupy a different world. 
 There seems to be a dependency on an evolved embodiment 
 which goes against the classical universal machine paradigm. 
 We develop these comments in more detail later in this section.
 
 For the mathematician, definability is the key to representation. 
 As previously mentioned, the language functions by representing  
 basic modes of using the informational content of the structure 
 over which the language is being interpreted. Very basic language 
 corresponds to classical computational relationships, and is 
 local in import. If we extend the language, for instance, by allowing 
 quantification, it still conveys information about an 
 algorithmic procedure for 
 accessing information.  The new element is that the information 
 accessed may now be emergent, spread across a range of regions 
 of the organism, its representation very much dependent on the 
 material embodiment, and with the information accessed 
 via finitary computational procedures which also depend 
 on the particular embodiment. One can observe this 
 preoccupation with the details of the embodiment in the 
 work of the neuro-scientist Antonio Damasio. 
 One sees this in the following 
 description from Damasio's book,  {\em 
 The Feeling Of What Happens}, of the kind of 
 mental recursions Steven Pinker 
 was referring to above \cite[p.170]{damasio99} : 
 \begin{quote}
 As the brain forms images of an object -- such as a face, a 
 melody, a toothache, the memory of an event -- and as the 
 images of the object affect the state of the organism, yet another 
 level of brain structure creates a swift nonverbal account of the 
 events that are taking place in the varied brain regions activated as 
 a consequence of the object-organism interaction. 
 The mapping of the object-related consequences occurs in first-order neural maps representing the 
 proto-self and object; the account of the causal relationship between object and organism can only be captured in second-order neural maps. \dots 
 one might say that the swift, second-order nonverbal account narrates 
 a story: that of the organism caught in the act of representing its own 
 changing state as it goes about representing something else.
 \end{quote}
 Here we see the pointers to the elements working against 
 the classical independence of the computational content 
 from its material host. We may have a mathematical precision 
 to the presentation of the process. But the presentation 
 of the basic information has to deal with emergence of a 
 possibly incomputable mathematical character, and 
 so has to be dependent on the material instantiation. 
 And the classical computation relative to such information, implicit in the 
 quotations from Pinker and Damasio, will need to work 
 relative to these material instantiations. The mathematics 
 sets up a precise and enabling filing system, telling the brain 
 how to work hierarchically through emergent informational 
 levels, within an architecture evolved over millions of years. 

 There is some recognition of this scenario in the current 
 interest in the evolution of hardware -- see, for example, 
 Hornby, Sekanina and Haddow
  \cite{HSH08}. We tend to  agree with Steven Rose 
  \cite{rose}:
  \begin{quote}
  Computers are designed, minds have evolved. 
  Deep Blue could beat Kasparov at a
game demanding cognitive strategies, but ask it to escape 
  from a predator, find food or a mate, and negotiate the 
  complex interactions of social life outside the chessboard or 
  express emotion when it lost a game, and it couldn't even leave 
  the launchpad. Yet these are the skills that human 
  survival depends on, the products of 3bn 
  years of trial-and-error evolution.
  \end{quote}
  From a computer scientist's perspective, we are grappling with 
the design of   
a {\em cyber-physical system} (CPS). And as 
Edward Lee from Berkeley describes \cite{lee08}:
\begin{quote}
 To realize the full potential of CPS, we will have to rebuild computing 
 and networking abstractions. These abstractions will have to 
 embrace physical dynamics and computation in a unified way.
 \end{quote}
 In Lee\cite{lee06}, he argues for ``a new systems science that is jointly physical and computational."
  
  Within such a context, connectionist 
  models with their close relationship to 
  synaptic interactions, and availability for ad hoc 
  experimentation, do seem to have a useful role. But their are good 
  reasons for looking for 
  a more fundamental mathematical model with which to 
  express the `design' on which to base a definable emergence. 
  The chief reason is the need for a general enough mathematical 
  framework, capable of housing different computationally complex 
  frameworks. Although the human brain is an important example, it is 
  but one part of a rich and heterogeneous computational 
  universe, reflecting in its workings many elements of that larger 
  context. The history of mathematics has led us to look for 
  abstractions which capture a range of related structures, 
  and which are capable of theoretical development informed 
  by intuitions from different sources,  which become 
  applicable in many different situations. And which provide 
  basic understanding to take us beyond the particularities 
  of individual examples.

 \section{Definability in What Structure?}
 
In looking for the mathematics to express the design, we need to take 
account of the needs of physics as well as those of mentality 
or biology. In his {\em The Trouble With Physics} \cite{smolin}, Lee Smolin 
points to a number of deficiencies of the standard model, and also 
of popular proposals such as those of string theory for filling  
its gaps. And in successfully modelling the physical universe, Smolin 
declares \cite[p.241]{smolin}:
\begin{quote}
{\em \dots causality itself is fundamental.}
\end{quote}
Referring to `early champions of the role of causality' such as 
Roger Penrose, Rafael Sorkin (the inventor of causal sets), Fay Dowker 
and Fotini Markopoulou, Smolin goes on to explain \cite[p.242]{smolin}:
\begin{quote}
It is not only the case that the spacetime geometry 
determines what the causal relations are. This can be turned around: 
Causal relations can determine the spacetime geometry \dots \par\noindent
It's easy 
to talk about space or spacetime emerging from something more fundamental, 
but those who have tried to develop the idea have found it 
difficult to realize in practice. \dots  We now believe they failed 
because they ignored the role that causality plays in spacetime. 
These days, many of us working on quantum gravity believe 
that causality itself is fundamental -- and is thus meaningful even 
at a level where the notion of space has disappeared.
\end{quote}
So, when we have translated `causality' into something meaningful, 
and the model based on it put in place --  the hoped-for 
prize is a theory in which even the background character of the universe is  
determined by its own basic structure. In such  a scenario, not only 
would one be able to do away with the need for exotic 
multiverse proposals, patched with inflationary theories and 
anthropic metaphysics. But , for instance, one can describe a proper 
basis for the variation of natural laws near a 
mathematical singularity, and so provide a mathematical foundation for 
the reinstatement of the philosophically more satisfying cyclical universe as an 
alternative to the inflationary big bang hypothesis  -- 
see Paul Steinhardt and Neil Turok's book \cite{ST07} for 
a well worked out proposal based on superstring theory. 
 
 \section{The Turing Landscape, Causality and Emergence \dots}
 
 If there is one field in which `causality' can be said to be fundamental, 
 it is that of computability. Although the sooner we can translate 
 the term into something more precise, the better. `Causality', despite 
 its everyday usefulness, on closer inspection is fraught with 
 difficulties, as John Earman \cite[p.5]{earman} nicely points out:
 \begin{quote}
\dots the most venerable of all the philosophical definitions 
[of determinism] holds that the world is deterministic 
just in case every event has a cause. The most immediate 
objection to this approach is that it seeks to explain a vague 
concept -- determinism -- in terms of a truly obscure one -- causation.
 \end{quote}
 Historically, one recognised the presence of a causal relationship 
 when a clear mechanical interaction was observed. But 
 Earman's book makes us aware of the subtleties beyond 
 this at all stages of history. The success of science in 
 revealing such interactions underlying mathematically 
 signalled   
 causality -- even for Newton's gravitational 
 `action at a distance' -- has encouraged us to 
 think in terms of mathematical relationships 
 being the essence of causality. Philosophically 
 problematic as this may be in general, there are enough 
 mathematical accompaniments to basic laws of nature 
 to enable us to extract a suitably general mathematical model of 
 physical causality. And to use this to improve our 
 understanding of more complicated (apparent) 
 causal relationships. The classical paradigm is 
 still Isaac Newton's formulation of a mathematically 
 complete formulation of his laws of motion, 
 sufficient to predict an impressive range 
 of planetary motions. 
 
 Schematically, logicians 
 at least have no problem representing Newtonian  
 transitions between mathematically well-defined 
 states of a pair of particles at different times as 
 the Turing reduction of one real to another, via 
 a partial computable (p.c.) functional describing what Newton said 
 would happen to the pair of particles. The functional 
 expresses the computational and continuous nature 
 of the transition. One can successfully use the functional 
 to approximate, to any degree of accuracy, a 
 particular transition. 
 
 This type of model, using partial computable functionals extracted 
 from Turing's \cite{Tu39} notion of oracle Turing machine, is 
 very generally applicable to basic laws of nature. 
 However, it is well-known that instances of a basic law can 
 be composed so as to get much more problematic mathematical 
 relationships, relationships which have a claim to be causal. 
 We have mentioned case above -- for instance those related to 
 the 3-body problem. Or strange attractors emergent from 
 complex confluences of applications of basic laws. 
 See recent work  Beggs, Costa, Loff and Tucker\cite{BCLT}, 
 Beggs and Tucker\cite{BT} 
 concerning the modelling of physical interactions as computation 
 relative to oracles, and 
 incomputability from mathematical thought experiments 
 based on Newtonian laws. 
 
 The technical 
 details of the extended Turing model, providing 
 a model of computable content of structures 
based on p.c.\ functionals over the reals, 
  can be found in Cooper\cite{Co04}. One can also find there 
  details of how Emil Post \cite{post48} used this model to 
  define the {\em degrees of unsolvability} -- now known as the 
  {\em Turing degrees} --  as a classification of 
  reals in terms of their relative computability. 
The resulting structure has turned out to be a very rich one, 
  with a high degree of structural pathology. 
 At a time when primarily mathematical motivations 
 dominated the field -- known for many years 
 as a branch of mathematical logic called 
 {\em recursive function theory} -- this pathology 
 was something of a disappointment. Subsequently, 
 as we see below, this pathology became the basis of a 
 powerful expressive language, delivering a 
 the sort of richness of definable relations 
 which qualify the structure for an important real-world 
 modelling role. 
 
 Dominant as this Turing model is, widely accepted to have a 
 canonical role, there are more general types of relative computation. 
 Classically, allowing non-deterministic Turing computations 
 relative to well-behaved oracles gives one nothing new. But 
 in the real world one often has to cope with data which is 
 imperfect, or provided in real-time, with delivery of 
 computations required in real time. There is an 
 argument that the corresponding generalisation is the 
 `real' relative computability. There are equivalent formalisations 
 -- in terms of {\em enumeration reducibility} between sets of data, 
 due to Friedberg and Rogers \cite{FR59}, or (see Myhill \cite{My61}), 
 in terms of {\em relative computability of partial functions} (extending 
 earlier notions of  Kleene and Davis). 
 The corresponding 
 extended structure provides an interesting and informative context 
 for the better known Turing degrees -- see, for example, Soskova and Cooper\cite{SC08}. 
 The Bulgarian research school, including D.\ Skordev, I.\ Soskov, 
 A.\ Soskova, A.\ Ditchev, H.\ Ganchev, M.\ Soskova and others has 
 played a special role in the development of the research area. 
 
 The universe we would like to model is one in which 
 we can describe global relations in terms of local structure --  
so capturing the emergence of large-scale formations, and 
 giving formal content to the intuition that such emergent 
 higher structures `supervene' on the computationally 
 more basic local relationships. 
 
Mathematically, there appears to be strong 
 explanatory power in the  formal modelling of 
 this scenario  
 as definability over a structure based on reducibilities closely 
 allied to 
  Turing functionals:  
or more generally, freeing the model from an explicit 
 dependence on language,  as Invariance under automorphisms 
 of the Turing structure. In the next section, we focus on the standard 
 Turing model, although the evidence is that similar outcomes 
 would be provided by the related models we have mentioned.

 \section{An Informational Universe, and Hartley Rogers' Programme}
 
 Back in 1967, the same year that Hartley Rogers' influential book 
 {\em  Theory of Recursive Functions and Effective Computability} appeared, 
 a paper \cite{Ro67}, based on an earlier talk of Rogers,  
 appeared in the proceedings volume of the 1965 Logic 
 Colloquium in Leicester.  This short article initiated a research agenda which 
 has held and increased its interest over a more than 40 year period. 
 Essentially, {\em Hartley Rogers' Programme} concerns the 
 fundamental problem of {\em characterising the Turing invariant relations}. 
 
 The intuition is that these invariant relations 
 are key to pinning down how basic laws and entities emerge as
mathematical constraints on causal structure. 
Where the richness of Turing structure discovered so far becomes the raw material for a multitude of non- trivially definable relations. There is an interesting 
 relationship here between the mathematics and the use of the anthropic 
 principle in physics to explain why the universe is as it is. It is well-known 
 that the development of the complex development we see around us is 
 dependent on  a subtle balance of natural laws and associated constants. 
 One would like the mathematics to explain why this balance 
 is  more than an accidental feature of one of a multitude, perhaps infinitely many, 
 randomly occurring universes. What the Turing universe delivers is a 
 rich infra-structure of invariant relations, providing a basis for a correspondingly 
 rich material instantiation, complete with emergent laws and constants,  
 a provision of strong determinism,  and a globally originating causality 
 equipped  with non-localism -- 
 though all in a very schematic framework. 
 Of course, echoing Smolin, it is  the underlying scheme that is currently missing. 
 We have a lot of detailed information, but the skeleton holding it all 
 together is absent. 
 
 However, the computability theorists have their own `skeleton in the cupboard'. 
 The modelling potential of the extended Turing model depends on it 
 giving some explanation of such well-established features as quantum 
 uncertainty, and certain experimentally verified uncertainties 
 relating to human mentality. And there is a widely believed mathematical 
 conjecture which would rob the Turing model of basic  
 credentials for modelling observable uncertainty. 
 
 The {\em Bi-Interpretability Conjecture}, arising from Leo Harrington's 
 familiarity with the model theoretic notion of 
 bi-interpretability, can be roughly described as asserting that:
 \begin{quote}
{\em  The Turing definable relations are exactly those with 
 information content describable in second-order arithmetic. }
 \end{quote}
 Moreover, given any description of  information content in second-order 
 arithmetic, one has a way of reading off the computability-theoretic definition in the 
 Turing universe. Actually, a full statement of the conjecture would be in terms 
 of `interpreting' one structure in another, a kind of poor-man's isomorphism. 
 Seminal work on formalising the global version of the conjecture, and 
 proving partial versions of it complete with key consequences and equivalences, 
 were due to Theodore Slaman and Hugh Woodin. See   
 Slaman's 1990 International Congress of Mathematicians 
 article \cite{slaman91} for a still-useful introduction to the conjecture 
 and its associated research project. 
 
 An unfortunate consequence of the conjecture being confirmed would be the 
 well-known rigidity of the structure second-order arithmetic being carried over to 
 the Turing universe. The breakdown of definability we see in the real world 
 would lose its model. However, work over the years makes this increasingly 
 unlikely. 
 
  See Nies, Shore and Slaman\cite{NSS98} for further development of   the requisite 
 coding techniques in the local context, with the establishment of a number 
 of local definability results. See Cooper\cite{Co99a,Cota} for work in the other 
 direction, both at the global and local levels. What is so promising here is 
 the likelihood of the final establishment of a subtle balance between 
 invariance and non-invariance, with the sort of non-trivial 
 automorphisms needed to deliver  a credible basis for the various 
 symmetries, and uncertainties peculiar to mentality and basic physics: along 
 with the provision via partial versions of bi-interpretability of an 
 appropriate model for the emergence of the more 
 reassuring `quasi-classical' world 
 from out of quantum uncertainty, and of other far-reaching consequences 
 bringing such philosophical concepts as epistemological relativism 
 under a better level of control. 
  
 {\bf To summarise:} What we propose is that this most cartesian of research areas, 
 classical computability theory, regain the real-world significance it 
 was born out of in the 1930s. And that it structure the informational 
 world of science in a radical and revealing way. The main features of this 
 informational world, and its modelling of the basic causal structure of the universe 
 would be:
 \begin{itemize}
 \item A universe described in terms of reals \dots
 \item With basic natural laws modelled by computable relations between reals. 
\item Emergence described in terms of definability/invariance over the resulting 
 structure \dots
 \item With failures of definable information content
modelling mental phenomena, quantum ambiguity, etc.  \dots
 \item Which gives rise to new levels of computable structure \dots 
 \item And a familiarly fragmented scientific enterprise. 
\end{itemize} 
 
 As an illustration of the explanatory power of the model, we 
 return to the problem of mental causation. 
 Here is William Hasker, writing in 
 {\em The Emergent Self} \cite[p.\ 175]{hasker99}, and trying to reconcile the 
 automomy of the different levels:
 \begin{quote}
The ``levels'' involved are levels of organisation and integration, and the downward influence means that the behavior of ``lower'' levels -- that is, of the components of which the ``higher-level'' structure consists -- is different than it would otherwise be, because of the influence of the new property that emerges in consequence of the higher-level organization.
\end{quote}
The mathematical model, making perfect sense of this, treats the brain and its 
emergent mentality as an organic whole. In so doing, it replaces the simple 
everyday picture of what a causal relationship is with a more subtle 
confluence of mathematical relationships. Within this confluence, one may for 
different purposes or necessities adopt different assessment of what the 
relevant causal relationships are. For us, thinking about this article, we 
regard the mentality hosting our thoughts to provide the significant 
causal structure. Though we know full well that all this mental activity 
is emergent from an autonomous brain, modelled with some validity via 
a neural network. 

So one might regard causality as a misleading concept in this context. 
Recognisable ÔcausalityÕ occurs at different levels of the
model, connected by relative definability.  And the causality 
at different levels in the form of relations with identifiable algorithmic
content, this content at higher levels being emergent. 
The diverse levels form a unity, with the ÔcausalÕ structure 
observed at one level reflected at other levels --  with the possibility of non-algorithmic ÔfeedbackÕ between levels. The 
incomputability involved in the transition between levels makes the 
supervenience involved have a non-reductive character.

\end{document}